\documentclass[12pt,a4paper,oneside]{amsart}
\usepackage[a4paper]{geometry}
\usepackage{mathpazo,perpage}
\newtheorem*{theorem}{Theorem}

\newcommand{\vertbar}{\>|\>}
\newcommand{\set}[2]{\ensuremath{\{ #1 \vertbar #2 \}}}

\DeclareMathOperator{\id}{id}

\MakePerPage[2]{footnote}
\begin{document}

\title{Melikyan algebra is a deformation of a Poisson algebra}

\author{Hayk Melikyan}
\address{Department of Mathematics and Physics, North Carolina Central University,
Durham, NC 27707, USA}
\email{melikyan@nccu.edu}

\author{Pasha Zusmanovich}
\address{Institute of Mathematics and Statistics, University of S\~{a}o Paulo,
S\~{a}o Paulo, Brazil}
\email{pasha.zusmanovich@gmail.com}
%\email{pasha@ime.usp.br}

\date{First written January 11, 2014; last minor revision April 24, 2016}
\thanks{J. Phys. Conf. Ser. \textbf{532} (2014), 012019; arXiv:1401.2566}

\begin{abstract}
We prove, using computer, that the restricted Melikyan algebra of dimension $125$
is a deformation of a Poisson algebra.
\end{abstract}

\maketitle

\section*{Introduction}

The classification of finite-dimensional simple Lie algebras over an 
algebraically close field of characteristic $p>3$ spanned the period of more than
50 years, and was achieved relatively recently 
(see \cite{premet-strade} and \cite{strade}).
Roughly, each simple algebra is either of classical type or of Cartan type.
As it is customary in many branches of mathematics, small characteristics
provide additional difficulties: in characteristic $p=5$, there is a series
of simple Lie algebras, discovered by the first author (\cite{melikyan}), and 
dubbed thereafter as Melikyan algebras (see \cite[Vol. I, \S 4.3]{strade}), 
which do not have analogs in higher characteristics. They are, however, 
intimately related to Lie algebras of Cartan type, and the present note provides
an additional evidence of such a connection. 

It was observed long time ago that simple Lie algebras of Cartan type can be 
deformed into each other (see, for example, \cite{dzhu} and references therein), and such kind of 
relationship, being interesting for its own sake, is also proved to be useful in
structure theory. A similar connection exists between Melikyan algebras and Lie
algebras of Cartan type: though they are different, Melikyan algebras are 
deformations of certain kind of Poisson algebras, the latter being 
extensions of Hamiltonian algebras. This was stated without proof in 
\cite[\S 5]{dk}. Unfortunately, the proof was never published, and more than 
decade later, the details seem to be lost.

In the present note we give a proof of this statement in the case of the smallest
(and also the only restricted one) algebra in the series of dimension $125$, 
$\mathcal M(1,1)$. The strategy of our proof is the following.
We construct on computer the deformed algebra $\mathcal D$ in question.
Then we verify, using computer, that the constructed Lie algebra
is simple. Over an algebraically closed field of characteristic $5$, there are 
only three simple Lie algebras of dimension $125$ -- 
the algebra of general type $\mathcal W_1(3)$, 
the contact algebra $\mathcal K_3(1,1,1)$, 
and the Melikyan algebra $\mathcal M(1,1)$, the latter has a symmetric
invariant bilinear form, and the formers do not. We check that the constructed
algebra $\mathcal D$ has a symmetric invariant bilinear form, and hence it is isomorphic
to the Melikyan algebra. 
After giving all necessary definitions and constructions in \S \ref{sec-main},
the computer calculations are explained in detail in \S \ref{sec-comp}.
In the last \S \ref{sec-next} we outline some further questions and possible lines of
development.

\section{Poisson algebras, cocycles and deformations}\label{sec-main}

We remind some standard definitions and facts concerning Poisson algebras, 
in a slightly more general, than usual, setting which will be useful for our 
purposes.

Let $A$ be a commutative associative algebra with unit $1$. 
Let us call a linear map $D: A \to A$ a \emph{generalized derivation} of $A$, 
if it satisfies the identity
$$
D(ab) = D(a)b + aD(b) - abD(1)
$$
for any $a,b \in A$. Note that this class of maps includes ordinary derivations 
(when $D(1) = 0$), and multiplications by an element of $A$. It is 
straightforward to check that the set of all generalized derivations of an
algebra $A$ forms a Lie algebra under the operation of commutator (which 
contains the usual derivations as a subalgebra).

Given a set $\{D_1, \dots, D_n, F_1, \dots F_n\}$ of pairwise commuting 
generalized derivations of $A$, the bracket
\begin{equation}\label{eq-poisson}
[a,b] = \sum_{i=1}^n D_i(a)F_i(b) - F_i(a)D_i(b)
\end{equation}
defines a Lie algebra structure on $A$ which will be denoted as 
$A_{\{D_i,F_i\}}$. In the particular case when
$A$ is a polynomial algebra in even number of variables 
$x_1, \dots, x_n, y_1, \dots, y_n$ (or, more generally, 
an algebra of smooth functions on a symplectic manifold), and generalized 
derivations are the usual ones, the partial derivatives with respect to 
$x_i, y_i$, $i=1, \dots, n$, the so obtained algebra is the ubiquitous Poisson 
algebra.

In fact, every summand in (\ref{eq-poisson}) defines a Lie algebra structure on
$A$. In other words, they form \emph{compatible Poisson structures}. This means
that first, for any two subsets 
$\set{D_i, F_i}{i\in I}$, $\set{D_i, F_i}{i\in J}$ of the initial set of 
commuting generalized derivations, the bracket defining the Lie algebra structure
on $A_{\set{D_i,F_i}{i\in J}}$, is a $2$-cocycle on the Lie algebra 
$A_{\set{D_i,F_i}{i\in I}}$, and second, the Massey square of this cocycle is 
zero (or, in other words, the infinitesimal deformation defined by this cocycle
is prolonged trivially).

Now, following \cite{dk}, consider a certain Poisson algebra structure on 
a $2$-variable divided powers algebra. Recall that the divided powers algebra
$\mathcal O_n(m_1, \dots, m_n)$, corresponding to $n$-tuple
of positive integers $m_1, \dots, m_n$, is a commutative associative
algebra defined over a field of characteristic $p>0$, with the basis 
$\set{x_1^{(i_1)} \cdots x_n^{(i_n)}}{0 \le i_k < p^{m_k}, k=1, \dots, n}$, 
and multiplication 
$$
x_1^{(i_1)} \cdots x_n^{(i_n)} \cdot x_1^{(j_1)} \cdots x_n^{(j_n)} = 
\binom{i_1+j_1}{i_1} \cdots \binom{i_n + j_n}{i_n} 
x_1^{(i_1+j_1)} \cdots x_n^{(i_n+j_n)} ,
$$
where the standing assumption is that $x_k^{(0)} = 1$, the unit of the algebra, 
and $x_k^{(i)} = 0$ if $i$ is outside the allowed range, i.e. 
$i<0$ or $i\ge p^{m_k}$, $k = 1, \dots, n$.

Obviously, 
\begin{equation}\label{eq-divpow}
\mathcal O_n(m_1, \dots, m_n) \simeq 
\mathcal O_1(m_1) \otimes \dots \otimes \mathcal O_1(m_n) .
\end{equation}

The algebra $\mathcal O_n(m_1, \dots, m_n)$ possess the following $n$ derivations, 
lowering the power of each variable: 
$$
\partial_{x_k} : 
x_1^{(i_1)} \cdots x_k^{(i_k)} \cdots x_n^{(i_n)} 
\mapsto  
x_1^{(i_1)} \cdots x_k^{(i_k-1)} \cdots x_n^{(i_n)} ,
$$
$k=1, \dots, n$. Obviously, these derivations pairwise commute. Of course,
any map of $\mathcal O_n(m_1, \dots, m_n)$ of the form $a\partial_{x_k}$, where
$a\in \mathcal O_n(m_1, \dots, m_n)$, is a derivation too (in fact, the 
derivation algebra of $\mathcal O_n(m_1, \dots, m_n)$ is linearly spanned by 
derivations of such form, and constitute the general Lie algebra of Cartan type 
$\mathcal W_n(m_1, \dots, m_n)$).

Let specialize the bracket (\ref{eq-poisson}) to the following situation:
$A = \mathcal O_2(2,1)$ (it will be more convenient to denote the 
generating variables in this case as $x$ and $y$ instead of $x_1$ and $x_2$),
$p=5$, $n=1$, and $D_1 = \partial_x$, $F_1 = \partial_y$. The resulting Poisson
Lie algebra $\mathcal P$ of dimension $125$ has one-dimensional center $K1$, 
and the commutator of the quotient $[\mathcal P/K1, \mathcal P/K1]$ is isomorphic to the simple 
$123$-dimensional Lie algebra $\mathcal H_2(2,1)$ of Hamiltonian type.

Now define the following $2$-cochains 
$\varphi, \psi: \mathcal P \times \mathcal P \to \mathcal P$:
\begin{align*}
\varphi(a,b) &= 
\partial_x^2 (a) \partial_x^3 (b) - \partial_x^3(a) \partial_x^2(b) \\
\psi(a,b) &= 
(\id - x \partial_x)(a) \partial_x^5 (b) - \partial_x^5(a) (\id - x\partial_x)(b)
\end{align*}
for $a,b\in \mathcal O_2(2,1)$. Here $\id$ denotes the identity map. Note that
both these cochains are of the form (\ref{eq-poisson}) for $n=1$, but while 
$\psi$ is formed by commuting generalized derivation $\id - x \partial_x$ and 
derivation $\partial_x^5$, and hence is automatically a cocycle on $\mathcal P$ 
by the discussion above, the map $\varphi$ is formed by the maps $\partial_x^2$ and 
$\partial_x^3$ which are not derivations. However, straightforward computations
(which may be performed on a computer, see below) show that $\varphi$ is a 
cocycle too.

Consider the cocycle $\varphi + 2\psi$. Direct computations (which, again, can
be performed with the aid of computer) show that this cocycle can be prolonged 
trivially, and hence define a deformation of the algebra $\mathcal P$. Let us 
denote this deformation by $\mathcal D$.

\begin{theorem}[Kostrikin--Dzhumadil'daev]
$\mathcal D \simeq \mathcal M(1,1)$.
\end{theorem}

This theorem is proved on computer, as explained in the next section.

\section{Computer calculations}\label{sec-comp}

We use GAP \cite{gap} which already contains a good deal of Lie-algebraic 
structures and functions for computation of various invariants of Lie algebras,
and augment them by the necessary stuff. All GAP routines described in this
section can be found at \newline
\texttt{http://justpasha.org/math/poisson-melikyan/}\footnote{
Added January 9, 2016: currently available as ancillary files accompanying the 
arXiv version of this paper.}.

The base field in Theorem is algebraically closed, of characteristic $5$. 
However, on computer we naturally work over the finite field $GF(5)$. 
While this difference is almost immaterial (all algebras in question are 
definable over $GF(5)$), in order to be rigorous, we distinguish between algebra
$\mathcal X$, where $\mathcal X \in \{\mathcal P, \mathcal D\}$, defined over 
an algebraically closed field $K$, and its
form $\mathcal X^\prime$ defined over the prime subfield $GF(5)$
(so $\mathcal X^\prime \otimes_{GF(5)} K \simeq \mathcal X$ as $K$-algebras).

\medskip

\emph{Step 1. Construction of the algebra $\mathcal P^\prime$}.

First, we provide two routines \texttt{PoissonAlgebra} and 
\texttt{SumOfAlgebraStructures} to construct Poisson algebras of kind 
(\ref{eq-poisson}). The routine \texttt{PoissonAlgebra} takes as arguments
an algebra $A$, and two linear maps $D,F: A \to A$, and returns the bracket
(\ref{eq-poisson}) in the particular case $n=1$: 
$$
[a,b] = D(a)F(b) - F(a)D(b)
$$
for $a,b\in A$. 
The routine \texttt{SumOfAlgebraStructures} takes as arguments several algebra 
structures defined on the same vector space $V$, and returns the algebra 
structure on $V$ which is the sum of the supplied algebra structures.

Second, we provide two routines \texttt{DividedPowersAlgebra} and 
\texttt{TensorProductOfAlgebras} to construct divided powers algebras.
The routine \texttt{DividedPowersAlgebra} takes as arguments a field $K$ of
positive characteristic, and a positive integer $n$, and returns the algebra 
$\mathcal O_1(n)$ over the field $K$. The routine 
\texttt{TensorProductOfAlgebras} takes as arguments several algebras defined 
over the same field, and returns their tensor product. Due to isomorphism (\ref{eq-divpow}), this enables us to construct
an arbitrary divided powers algebra.

All this, together with a few other auxiliary routines implementing the 
derivations $\partial_x$ and $\partial_y$ as linear operators acting on 
$\mathcal O_2(2,1)$, enables us to construct the algebra $\mathcal P^\prime$.

\medskip

\emph{Step 2. Construction of the algebra $\mathcal D^\prime$}.

Using the already employed routines \texttt{PoissonAlgebra} and 
\texttt{SumOfAlgebraStructures}, we construct cocycles $\varphi$ and $\psi$, 
and then the bracket
$$
\{a,b\} = [a,b] + \varphi(a,b) + 2\psi(a,b) ,
$$
where $a,b\in \mathcal O_2(2,1)$, what gives us the algebra $\mathcal D^\prime$.

\medskip

\emph{Step 3. Identifying algebras $\mathcal D$ and $\mathcal M(1,1)$}.

Generally, to establish on computer whether two given Lie algebras are 
isomorphic or not, amounts to solving a system of quadratic equations, and as 
such, is a difficult task. There are several algorithms for doing that, some of 
them quite sophisticated (see, for example, \cite{eick}), but neither of them,
as of time of this writing, would work for algebras of dimension 
$125$ on a reasonable computer within a reasonable amount of time.

However, to establish whether two modules over the \emph{same} algebra are 
isomorphic or not, is a linear problem, and as such, is much more tractable.
Some little theory enables us to take advantage of the latter fact.

In dealing with modules over algebras, we use a highly efficient Meataxe suit
of algorithms (see \cite[\S 7.4]{handbook}).
Initially, Meataxe was developed for dealing with
finite groups representations, but it works on the level of associative matrix
algebras (via group algebras), so it can be utilized for study of 
representations of Lie algebras (via associative envelopes of Lie algebras) 
equally well.
Using Meataxe routines available in GAP, we can compute whether a given Lie
algebra is simple (i.e., the adjoint representation is irreducible), 
central simple (i.e., the adjoint representation is absolutely irreducible), and 
possess a symmetric invariant bilinear form (i.e., the adjoint representation is 
equivalent to its dual). We borrowed the idea to use Meataxe in the 
Lie-algebraic context from \cite{eick}.

First, we establish that the algebra $\mathcal D^\prime$ is central simple
(so $\mathcal D$ is simple).
According to the classification of simple modular Lie algebras, over an 
algebraically closed field of characteristic $5$ there are only three simple Lie 
algebras of dimension $125$ -- the algebra of general type $\mathcal W_1(3)$,
the algebra of contact type $\mathcal K_3(1,1,1)$, and the Melikyan algebra 
$\mathcal M(1,1)$ (see, for example, \cite[Vol. I, \S 4.2]{strade} for the
dimensions of Lie algebras of Cartan type). In order to distinguish between
them, we have to choose some invariant, and the suitable invariant in our case 
is the presence of a nonzero symmetric invariant bilinear form. 

Recall that a symmetric bilinear form $\omega: L \times L \to K$ on a Lie 
algebra $L$ is called \emph{invariant}, if 
$\omega([x,z],y) + \omega(x,[y,z]) = 0$ for any $x,y,z\in L$. On simple Lie 
algebras, this form, if exists, is unique up to a scalar due to Schur lemma, and
the existence of nonzero form is equivalent to isomorphism of $L$-modules 
$L \simeq L^*$. The algebras $\mathcal W_1(3)$ and $\mathcal K_3(1,1,1)$ do not 
have such a form (see \cite[\S 4.6, Theorems 6.3 and 6.6]{sf}, and also 
\cite{invar-forms} for a short alternative proof and further references), while 
the Melikyan algebra does (\cite[Proposition 6.1]{premet-strade}).

As the last step, we establish that $\mathcal D^\prime$ has a symmetric 
invariant bilinear form, and hence so does $\mathcal D$. To summarize:
$\mathcal D$ is a simple Lie algebra of dimension $125$ having nonzero 
symmetric invariant bilinear form, and hence have to be isomorphic to the 
Melikyan algebra, thus proving the Theorem.

The whole computation takes about 1.5 minutes on a machine with 2.40 GHz CPU.

\section{What next?}\label{sec-next}

It goes without saying that a more satisfactory proof of the Theorem would cover the whole series
of Melikyan algebras, will not involve computer, and will provide an explicit
isomorphism. We are planning to systematically compute deformations of Poisson 
Lie algebras of the kind appearing in this note, and identify among them those
which are isomorphic to Melikyan algebras and other interesting Lie algebras in
small characteristics. As a first step, we envisage
a linear-algebraic approach to computation of (low-degree) cohomology of a more 
general class of Poisson algebras, in which the underlying commutative 
associative algebra $A$ is decomposed as the tensor product of two algebras. This
should be somewhat similar to computation of cohomology of current Lie 
algebras in \cite{low}.

Another interesting question is how the $p$-map comes into play. Note that the 
Poisson algebra $\mathcal P$ is not restricted, while its deformation 
$\mathcal M(1,1)$ is. This phenomenon, where a restricted Lie algebra is a 
deformation of a not restricted one, and, generally, how the $p$-map behaves
under deformations, deserves further investigation.

\section*{Acknowledgements}

Thanks are due to Askar Dzhumadil'daev for turning our attention to \cite{dk}
in the first place, and for subsequent useful discussion. This work was 
carried out while Zusmanovich was visiting North Carolina Central University, 
and was supported by the grants HRD-0833184 (NSF), NNX09AV07A (NASA), and 
ETF9038 (Estonian Science Foundation).

\end{document}